\documentclass[12pt]{article}
\title{Linnik's approximation to Goldbach's conjecture, and other problems}
\author{D.J. Platt\\ Heilbronn Institute for Mathematical Research \\ University of Bristol, Bristol, UK\\ dave.platt@bris.ac.uk\\ \\
and\\ \\
T.S. Trudgian\footnote{Supported by Australian Research Council DECRA Grant DE120100173.}\\
Mathematical Sciences Institute\\ The Australian National University,
 ACT 0200, Australia\\ timothy.trudgian@anu.edu.au
}

\usepackage{comment}
\usepackage{url}
\usepackage{amsthm}
\usepackage{fullpage}
\usepackage{amsmath}
\usepackage{amssymb}
\usepackage[boxed]{algorithm2e}

\begin{document}
\maketitle

\begin{abstract}
\noindent
We examine the problem of writing every sufficiently large even number as the sum of two primes and at most $K$ powers of 2. We outline an approach that only just falls short of improving the current bounds on $K$. Finally, we improve the estimates in other Waring--Goldbach problems.
\end{abstract}
\section{Introduction}
The Goldbach conjecture is that every even $n>4$ can be written as a sum of two prime numbers. Linnik proved that there exists a finite $K$ such that, for all sufficiently large even $n$, one may write
\begin{equation}\label{link}
 n = p + q + 2^{\nu_{1}} + 2^{\nu_{2}} + \cdots + 2^{\nu_{r}},
 \end{equation}
where $p$ and $q$ are primes, the $\nu_{i}$ are positive integers, and where $r\leq K$. For a historical development on bounding the size of $K$, see \cite[\S 1]{HBP}.

Significant improvements on the size of $K$ were made by Heath-Brown and Puchta \cite{HBP} and, independently, by Pintz and Ruzsa \cite{PintzR}. Heath-Brown and Puchta showed that $K\leq 13$, and on the Generalised Riemann Hypothesis (GRH) that $K\leq 7$. Pintz and Ruzsa established this latter result and announced that they could show unconditionally that $K\leq 8$. This paper is yet to appear in print.
Elsholtz, in an unpublished manuscript, showed that $K\leq 12$; this was proved independently by Liu and L\"{u} \cite{LiuL}.

The methods of Heath-Brown and Puchta and of Pintz and Ruzsa allow one to show that $K$ is admissible in (\ref{link}) provided that
\begin{equation}\label{12.98}
\lambda^{K-2} < \frac{C_{3}}{(C_{1} - 2)C_{2} + c C_{0}^{-1}\log 2},
\end{equation}
for certain constants $C_{0}, C_{1}, C_{2}, C_{3}, \lambda$ and $c$. The inequality (\ref{12.98}) follows from \S 6 in \cite{HBP} and the relations (8.20) and (8.21) in \cite{PintzR}, which, although proved under GRH, give analogous results unconditionally.

We have
\begin{equation*}\label{C0}
C_{0} = \prod_{p>2} (1- (p-1)^{-2}),
\end{equation*}
where, according to Wrench \cite{Wrench}, 
\begin{equation}\label{C0bound}
0.6601618158< C_{0} < 0.6601618159.
\end{equation}
Indeed, Wrench computed $C_{0}$ to 45 digits; the truncated bound in (\ref{C0bound}) is certainly fit for purpose.

As for $C_{2}$ we have
\begin{equation*}\label{C2}
C_{2} = \sum_{d=1}^{\infty} \frac{k(2d-1)}{\epsilon(2d-1)},
\end{equation*}
where $k(n)$ is a multiplicative function defined by
\begin{equation*}
k(p^{e}) =\begin{cases} 0, \quad &p=2\; \textrm{or } e\geq 2,\\
(p-2)^{-1}, \quad &\textrm{otherwise},
\end{cases}
\end{equation*}
and where $\epsilon(d)$ is the multiplicative order of $2$ modulo $d$. Lemma $2''$ of \cite{PintzRom} shows that $1.2783521041< C_{2}C_{0}< 1.2784421041$, which may be combined with (\ref{C0bound}) to give
\begin{equation*}\label{C2bound}
1.93642< C_{2} < 1.93656.
\end{equation*}

The constant $C_{1}$ is that which appears in Chen's work on the twin prime conjecture. Heath-Brown and Puchta take $C_{1} = 7.8342$ as in \cite{Chen}; Liu and L\"{u} take $C_{1} = 7.8209$ as in \cite{Wu}.\footnote{It is claimed in \cite{Donghua} that one may take $C_{1} = 7.81565$; this proof appears to be incomplete.} 

The constant $C_{3}$ satisfies
\begin{equation*}\label{C3}
2 \sum_{d\leq D} k(d) H(d;N, K) \epsilon(d)^{-K} \geq C_{3},
\end{equation*}
in which $D$ is a fixed number, and 
$$H(d; N, K) = \#\{(\nu_{1}, \ldots, \nu_{K}): 1\leq \nu_{i}\leq \epsilon(d), d|N - \sum 2^{\nu_{i}}\}.$$
Heath-Brown and Puchta used $D=5$ to show that $C_{3} \geq 2.7895$; Liu and L\"{u} used $D=11$ to show that $C_{3} \geq 2.8096$; it is remarked in \cite{HBP} that Elsholtz used $D=21$ to show that $C_{3} \geq 2.96169.$

To define $c$ and $\lambda$ we first define
\begin{equation}\label{joint}
G_{L}(x) = \sum_{j=0}^{L-1} e(2^{j} x), \quad L = \lfloor \log N/\log 2 \rfloor,
\end{equation}
where $e(x) = \exp(2\pi i x)$. For $c>0$ we seek to find the smallest positive value of $\lambda$ such that
\begin{equation}\label{delta}
\Delta = \textrm{meas}(x \in[0,1]:\, |G_{L}(x)| > \lambda L) < N^{-\frac{c}{\log 2}}.
\end{equation}
Heath-Brown and Puchta showed that one may take $c = 109/154$ unconditionally; Pintz and Ruzsa claim that one may take $c= 3/5$. Both sets of authors show that one may take $c = \frac{1}{2}$ on GRH.
Heath-Brown and Puchta showed that $\lambda \leq 0.863665$; Liu and L\"{u} improved this to $\lambda \leq 0.862327.$ On GRH, Heath-Brown and Puchta showed that $\lambda \leq 0.722428$. Pintz and Ruzsa used an improved method to show that $\lambda \leq 0.716344.$

Using 
$$(C_{0}, C_{1}, C_{2}, C_{3}, \lambda) = (0.6601618159, 7.8209, 1.93656, 2.96169, 0.862327),$$
and the same quintuple with the last entry replaced by 0.716344, one may see that (\ref{12.98}) is satisfied for $K\geq 11.4549$ and $K \geq 6.1432$. This proves only what is already known, viz.\ that one may take $K= 12$, and $K= 7$ on GRH. 

It is clear that there is little to be gained by pursuing improvements in $C_{0}$ and in $C_{2}$. It appears difficult to improve $C_{1}$ at all substantially --- see \cite[Rem.\ 2, p.\ 253]{Wu}. In \S \ref{beef} we improve the value of $C_{3}$; in \S \ref{pork} we investigate $\lambda$. Finally, in \S \ref{other} we improve on estimates for some related problems.

\section{Computing $C_{3}$}\label{beef}
Heath-Brown and Puchta examined all those $d\leq D$ for which $2$ is a primitive root modulo $d$. They stated that, in this case,
\begin{equation*}\label{HBPs}
H(d; N, K) = \begin{cases} \frac{1}{d}\{(d-1)^{K} - (-1)^{K}\}, \quad & d \nmid N,\\
 \frac{1}{d}\{(d-1)^{K} + (-1)^{K}(d-1)\}, \quad & d \mid N.
 \end{cases}
 \end{equation*}

To improve on the value we could take for $C_3$ we wished to consider more general $d$. Algorithm \ref{alg:1} describes the approach we adopted.

\begin{algorithm}\label{alg:1}
\SetAlgoLined
\SetKwInOut{Input}{input}
\Input{$d>1$, an odd integer, and $K>1$ an integer}
$W\leftarrow \text{a vector of length}\,d\, \text{indexed by}\, W[0]\ldots W[d-1]$\;
$Y\leftarrow \text{a vector of length}\,d\, \text{indexed by}\, Y[0]\ldots Y[d-1]$\;
\For{$i\leftarrow 0$ \KwTo $d-1$}{
$W[i]\leftarrow \#\{\nu|1\leq \nu \leq\epsilon(d),2(K-1)+2^\nu \equiv i\,(\textrm{mod}\, d)\}$}
\For{$k\leftarrow 2$ \KwTo $K$}{
\For{$j\leftarrow 2$ \KwTo $\epsilon(d)-1$}{
$r\leftarrow 2^j-2$\;
\For{$i\leftarrow 0$ \KwTo $d-1$}{
$Y[i+r\,(\textrm{mod}\,d)]\leftarrow W[i]$}
\For{$i\leftarrow 0$ \KwTo $d-1$}{
$W[i]\leftarrow W[i]+Y[i]$}}}
$res\leftarrow W[0]$\;
\For{$i\leftarrow 1$ \KwTo $d-1$}{
\If{$W[i]<res$}{$res\leftarrow W[i]$}}
\textbf{return} $res$\;
\caption{Computing $H(d;N,K)$.}
\end{algorithm}

Setting up the initial vector requires $d+\epsilon(d)=\mathcal{O}(d)$ steps. We then copy and add the $d$ vector entries $\epsilon(d)$ times for each of the remaining $K-1$ powers of $2$ giving a total cost of $\mathcal{O}(Kd^2)$ for each $d$. Computing $H$ for all the $d$ less than $D$ therefore has time complexity of $\mathcal{O}(KD^3)$ where each operation is on a number of size $\mathcal{O}(D\log k)$ bits. The space requirement is $\mathcal{O}(D^2\log k)$.

Following a suggestion by Roger Heath-Brown, we introduce a slight variation of the above argument. Let the \textit{worst} residue class modulo $d$ be that which contributes the least to $H(d; N, K)$. Suppose, for example, that the worst residue class modulo 15 is $N\equiv 0$, and the worst residue class modulo 3 is $N \equiv 1$. Since there can be no values of $N$ that belong to both residue classes, this `worst of the worst' scenario does not arise. We therefore limit ourselves to \textit{admissible} values: that is, sets of residue classes for $N$ that  could be simultaneously satisfied.
 Modifying the above Algorithm \ref{alg:1} to take advantage of this is a trivial matter.

\subsection{Implementation and Results}

We implemented the above algorithm in C++ using GMP \cite{GMP2013} for the large integer arithmetic and MPFI \cite{Revol2005} to handle floating point quantities as intervals, thereby avoiding any issues with rounding. We summed over all $d\leq 40,000$ and considered admissibility modulo $255,255$ to obtain $C_3\geq 3.011112$ for $K=6$ and $C_3\geq 3.02858417$ for $K=11$. The computations required $58$ hours in the case of $K=6$ and $116$ hours in the case of $K=11$ on a single core of a $1.8$ GHz Intel\textsuperscript{\circledR} Xeon\textsuperscript{\circledR} E5-2603.

By modern standards these run times are modest; it would be a simple matter to go higher in $d$. However, we expect the returns to be very small based on the following argument. We shall ignore admissibility since the improvement we observed when introducing it to our algorithm was small. In this case, the best
 we can hope for is an even distribution of the counts over all the residue classes modulo $d$. Thus we expect $H(d;N,K)$ to be about $\epsilon(d)^{K}/d$, whence the contribution from each $d$ to $C_3$ will be no more than $\frac{2k(d)}{d}.$
We note that this treatment removes the dependency on $K$. Further, we have trivially that
$k(d)\leq (d-2)^{-1}$
for all $d>3$. Therefore
\begin{equation*}
2\sum\limits_{d>40,000}\frac{k(d)}{d}<2\int\limits_{40,000}^\infty \frac{dt}{t(t-2)}< 5.1\times 10^{-5}.
\end{equation*}
It seems that the potential gains from further computation are limited.

\section{Computing $\lambda$ and $K$}\label{pork}
To estimate $\lambda$, we use the method given by Pintz and Ruzsa \cite{PintzR}. Essentially one wishes to approximate $G_{L}(x)$ in (\ref{joint}) by estimating the error in $G_{2^{h}}(x) - G_{L}(x)$ for $h$ large. It is to this purpose that \S\S 3-7 of \cite{PintzR} are dedicated. Heath-Brown and Puchta took $h= 16$; Liu and L\"{u} took $h= 23$. Using the method of Pintz and Ruzsa we are able to take $h$ considerably larger. 

The authors are grateful to Alessandro Languasco who supplied his Pari \cite{pari} implementation of this algorithm. We took $h=2^{138}$ with polynomials\footnote{For details of these polynomials the reader is invited to examine \S 6 in \cite{PintzR}.}  of degree $40$ to obtain $\lambda\leq 0.8594000$ unconditionally and $\lambda\leq 0.7163436$ on GRH. The computations take about $15$ minutes using Pari, or about double that when implemented in C using the interval arithmetic package ARB \cite{johansson2014}. The latter approach confirms that the stated values for $\lambda$ are accurate to the precision given. 

With $(C_{3}, \lambda) = (3.02858417, 0.8594000)$ in the unconditional case and with $(C_{3}, \lambda)=(3.011112, 0.7163436)$ on GRH we have
\begin{equation}\label{close}
K \geq 11.0953 \quad \textrm{unconditionally}, \quad\quad\quad K\geq 6.09353 \quad\textrm{on GRH}.
\end{equation}
This means we are just short of being able to take $K=11$ and, on GRH, $K=6$. Given the difficulty in improving the values of $C_{3}$ and $\lambda$ with existing methods, it seems that a new idea is needed to improve the estimate on $K$.

As a consolation prize, we applied the same code to some other problems. For these problems one has a possibly different value of $c$ for which one wishes to calculate a small value of $\lambda$. The results are summarised in Table \ref{table3}.

\section{Other Waring--Goldbach problems}\label{other}
Suppose it is conjectured that for all sufficiently large $N$ we have $N= f_{1} + \cdots + f_{r}$ for certain numbers $f_{i}$. Suppose that we can prove the following approximation of this conjecture, that $N = f_{1} + \cdots + f_{r} + 2^{\nu_{1}} + \cdots + 2^{\nu_{r}}$, where $r\leq \mathcal{K}$ for some $\mathcal{K}$. Just as in Linnik's approximation to Goldbach's conjecture, one seeks good bounds on $\mathcal{K}$. Various approximations have been given to problems involving sums of powers of primes. We investigate ten of them below.

For the following problems the value of $c$ has been established in the literature. It may be possible to improve this value and some of the other arguments that lead to the estimates on $\mathcal{K}$ in problems (A)-(J). We have not pursued this: we limited ourselves to improving the value of $\lambda$ since this appears to be the most influential parameter.
\subsubsection{Even numbers as sums of four squares of primes} 
\begin{equation}\tag{A}
N = p_{1}^{2} + p_{2}^{2} + p_{3}^{2} + p_{4}^{2} + 2^{\nu_{1}} + 2^{\nu_{2}}+ \cdots + 2^{\nu_{K_{A}}}.
\end{equation}
This was considered in \cite{LLZ99, LiuLiu00, LiuLu04, Li06} and most recently by Zhao \cite{Zhao13} who showed that $K_{A}\leq 46.$
\subsubsection{Odd numbers as sums of a prime and two squares of primes}
\begin{equation}\tag{B}
N= p_{1} + p_{2}^{2} + p_{3}^{2} + 2^{\nu_{1}} + 2^{\nu_{2}}+ \cdots + 2^{\nu_{K_{B}}}.
\end{equation}
This was considered in \cite{LLZ99, Liu04, HuYang11, Li07, LuSun09} and most recently by Liu \cite{LiuB:2014} who showed that $K_{B} \leq 35$.
\subsubsection{Even numbers as sums of eight cubes of primes}
\begin{equation}\tag{C}
N = p_{1}^{3} + \cdots + p_{8}^{3} + 2^{\nu_{1}} + 2^{\nu_{2}}+ \cdots + 2^{\nu_{K_{C}}}.
\end{equation}
This was considered in \cite{LiuLu10} and most recently by Liu \cite{LiuDensity12} who showed that $K_{C}\leq 341.$
\subsubsection{Odd numbers as sums of a prime and four cubes of primes}
\begin{equation}\tag{D}
N= p_{1} + p_{2}^{3} + \cdots + p_{5}^{3} + 2^{\nu_{1}} + 2^{\nu_{2}} +\cdots + 2^{\nu_{K_{D}}}.
\end{equation}
This considered by Liu and L\"{u} \cite{LiuLu112} who showed that $K_{D}\leq 106$.

\subsubsection{Even numbers as sums of two squares of primes and four cubes of primes}
\begin{equation}\tag{E}
N= p_{1}^{2} + p_{2}^{2} + p_{3}^{3} \cdots + p_{6}^{3} + 2^{\nu_{1}} + 2^{\nu_{2}}+ \cdots + 2^{\nu_{K_{E}}}.
\end{equation}
This was considered by Liu and L\"{u} \cite{LiuLu112} who showed that  $K_{E} \leq 211$.
\subsubsection{Even numbers as sums of a prime, a square of a prime and two cubes of primes}
\begin{equation}\tag{F}
N= p_{1} + p_{2}^{2} + p_{3}^{3} + p_{4}^{3} + 2^{\nu_{1}} + 2^{\nu_{2}}+ \cdots + 2^{\nu_{K_{F}}}.
\end{equation}
This was considered by Liu and L\"{u} \cite{LiuLu11} who showed that $K_{F}\leq 161$.

In the nextfour problems one asks when the equations are true simultaneously for positive even integers $B_{1}$ and $B_{2}$ with $B_{1}>B_{2}$.
\subsubsection{Even numbers as sums of two primes, simultaneously}
\begin{equation}\tag{G}
\begin{split}
B_{1} &= p_{1} + p_{2} + 2^{\nu_{1}} + 2^{\nu_{2}} +\cdots + 2^{\nu_{K_{G}}}\\
B_{2} &= p_{3} + p_{4} + 2^{\nu_{1}} + 2^{\nu_{2}} +\cdots + 2^{\nu_{K_{G}}}.
\end{split}
\end{equation}
This was considered by Kong \cite{Kong13} who showed that $K_{G} \leq 63$ and that, on GRH, $K_{G} \leq 31$.
\subsubsection{Even numbers as sums of four squares of primes, simultaneously}
\begin{equation}\tag{H}
\begin{split}
B_{1} &= p_{1}^{2} + p_{2}^{2} + p_{3}^{2} + p_{4}^{2} + 2^{\nu_{1}} + 2^{\nu_{2}}+ \cdots + 2^{\nu_{K_{H}}}\\
B_{2} &=  p_{5}^{2} + p_{6}^{2} + p_{7}^{2} + p_{8}^{2} + 2^{\nu_{1}} + 2^{\nu_{2}}+ \cdots + 2^{\nu_{K_{H}}}.
\end{split}
\end{equation}
This was considered in \cite{Liu13} and most recently by Hu and Liu  \cite{HuLiu:2014} who showed that $K_{H}\leq 142$.\footnote{Though it seems that Hu and Liu's proof actually gives $K_{H} \leq 141$.}
\subsubsection{Even numbers as sums of eight cubes of primes, simultaneously}
\begin{equation}\tag{I}
\begin{split}
B_{1} &= p_{1}^{3} +\cdots + p_{8}^{3} + 2^{\nu_{1}} + 2^{\nu_{2}}+ \cdots + 2^{\nu_{K_{I}}}\\
B_{2} &=  p_{9}^{3} +\cdots + p_{16}^{3} + 2^{\nu_{1}} + 2^{\nu_{2}}+ \cdots + 2^{\nu_{K_{I}}}.
\end{split}
\end{equation}
This was considered by Liu \cite{Liu13} who showed that $K_{I}\leq 1432$.\footnote{Note that, on \cite[p.\ 3347]{Liu13} Liu uses $b= 268096$, which comes from Ren \cite{Ren03}. This value has been improved in \cite{LiuDensity12} to $b = 147185.22$. This gives at once that $K_{I} \leq 1364$.}

\subsubsection{Odd numbers as sums of one prime and two squares of primes, simultaneously}
\begin{equation}\tag{J}
\begin{split}
B_{1} &= p_{1} + p_{2}^{2} + p_{3}^{2} + 2^{\nu_{1}} + 2^{\nu_{2}}+ \cdots + 2^{\nu_{K_{J}}}\\
B_{2} &=  p_{4} + p_{5}^{2} + p_{6}^{2} +2^{\nu_{1}} + 2^{\nu_{2}}+ \cdots + 2^{\nu_{K_{J}}}.
\end{split}
\end{equation}
This was considered by Liu \cite{Liu:2013J} who showed that $K_{J}\leq 332$.

\begin{table}[ht] 
\caption{Improvements on (A)-(J)} 
\label{table3} 
\centering 
\begin{tabular}{c c c c c c} 
\hline\hline 
  & Required $c$ & Old $\lambda$ & New $\lambda$ & Old $\mathcal{K}$ & New $\mathcal{K}$ \\[0.5 ex] \hline 
(A) & 3/4 & 0.887167 & 0.8844473 & 46 &45\\ 
(B) & 3/4 & 0.887167   & 0.8844473 & 35 &34\\ 
(C) & 19/21 & 0.965411   & 0.9642399 & 341 &330\\ 
(D) & 113/126 & 0.961917 & 0.9606646 & 106 &102\\ 
(E) &53/63 & 0.935746    & 0.9339489 & 211 &205\\ 
(F) &109/126 & 0.947313  & 0.9457435 & 161&156\\ 
(G) &109/154 & 0.862327      & 0.8594000 & 63 &62\\ 
 (G) on GRH & 1/2 & 0.716344 & 0.7163436 & 31 & 31 \\ 
(H) & 3/4 & 0.887167  & 0.8844473 & 142 &138\\ 
(I) &19/21 & 0.965411 & 0.9642399 & 1432 &1319\\ 
(J) & 3/4 & 0.887167  & 0.8844473 & 332 &323\\ 
\hline\hline 
\end{tabular} 
\end{table} 

\bibliographystyle{plain}
\bibliography{themastercanada}

\begin{thebibliography}{10}

\bibitem{pari}
C.~Batut, K.~Belabas, D.~Bernardi, H.~Cohen, and M.~Olivier.
\newblock {User's Guide to PARI-GP}, 2000.

\bibitem{Chen}
J.-R. Chen.
\newblock On the {G}oldbach's problem and the sieve methods.
\newblock {\em Sci. Sinica}, 21(6):701--739, 1978.

\bibitem{GMP2013}
Torbj\"orn Granlund.
\newblock {\em {The GNU Multiple Precision Arithmetic Library}}, 5.1.1 edition,
  February 2013.

\bibitem{HBP}
D.~R. Heath-Brown and J.-C. Puchta.
\newblock Integers represented as a sum of primes and powers of two.
\newblock {\em Asian J. Math.}, 6(3):535--565, 2002.

\bibitem{HuLiu:2014}
L.~Hu and H.~Liu.
\newblock On pairs of four prime squares and powers of two.
\newblock {\em J. Number Theory}, 147:594--604, 2015.

\bibitem{HuYang11}
L.~Hu and L.~Yang.
\newblock The number of powers of 2 in a representation of large odd integers.
\newblock {\em Acta Arith.}, 150(2):175--192, 2011.

\bibitem{johansson2014}
Fredrik Johansson.
\newblock {Arb: a C library for ball arithmetic}.
\newblock {\em ACM Communications in Computer Algebra}, 47(3/4):166--169, 2014.

\bibitem{Wrench}
J.~W.~Wrench Jr.
\newblock Evaluation of {A}rtin's constant and the twin-prime constant.
\newblock {\em Math. Comp.}, 15(76):396--398, 1961.

\bibitem{Kong13}
Y.~Kong.
\newblock On pairs of linear equations in four prime variables and powers of
  two.
\newblock {\em Bull. Aust. Math. Soc.}, 87:55--67, 2013.

\bibitem{Li06}
H.~Z. Li.
\newblock Four prime squares and powers of 2.
\newblock {\em Acta Arith.}, 125:383--391, 2006.

\bibitem{Li07}
H.~Z. Li.
\newblock Representation of odd integers as the sum of one prime, two squares
  of primes and powers of 2.
\newblock {\em Acta Arith.}, 128:223--233, 2007.

\bibitem{LiuLu04}
J.~Liu and G.~L\"{u}.
\newblock Four squares of primes and 165 powers of 2.
\newblock {\em Acta Arith.}, 114(1):55--70, 2004.

\bibitem{LiuLiu00}
J.~Y. Liu and M.~C. Liu.
\newblock Representation of even integers as sums of squares of primes and
  powers of 2.
\newblock {\em J. Number Theory}, 83:202--225, 2000.

\bibitem{LLZ99}
J.~Y. Liu, M.~C. Liu, and T.~Zhan.
\newblock Squares of primes and powers of 2.
\newblock {\em Monatsh. Math.}, 128:283--313, 1999.

\bibitem{Liu04}
T.~Liu.
\newblock Representation of odd integers as the sum of one prime, two squares
  of primes and powers of 2.
\newblock {\em Acta Arith.}, 115:97--118, 2004.

\bibitem{LiuDensity12}
Z.~Liu.
\newblock Density of the sums of four cubes of primes.
\newblock {\em J. Number Theory}, 132:735--747, 2012.

\bibitem{Liu:2013J}
Z.~Liu.
\newblock On pairs of one prime, two prime squares and powers of $2$.
\newblock {\em Int. J. Number Theory}, 9(6):1413--1421, 2013.

\bibitem{Liu13}
Z.~Liu.
\newblock On pairs of quadratic equations in primes and powers of 2.
\newblock {\em J. Number Theory}, 133:3339--3347, 2013.

\bibitem{LiuB:2014}
Z.~Liu.
\newblock One prime, two squares of primes and powers of $2$.
\newblock {\em Acta Math. Hungar.}, 143(1):3--12, 2014.

\bibitem{LiuLu10}
Z.~Liu and G.~L\"{u}.
\newblock Eight cubes of primes and powers of 2.
\newblock {\em Acta Arith.}, 145(2):171--192, 2010.

\bibitem{LiuL}
Z.~Liu and G.~L\"{u}.
\newblock Density of two squares of primes and powers of 2.
\newblock {\em Int. J. Number Theory}, 7(5):1317--1329, 2011.

\bibitem{LiuLu112}
Z.~Liu and G.~L\"{u}.
\newblock Two result on powers of 2 in {W}aring--{G}oldbach problem.
\newblock {\em J. Number Theory}, 131:716--736, 2011.

\bibitem{LiuLu11}
Z.~X. Liu and G.~L\"{u}.
\newblock On unlike powers of primes and powers of 2.
\newblock {\em Acta Math. Hungar.}, 132(1-2):125--139, 2011.

\bibitem{LuSun09}
G.~L\"{u} and H.~Sun.
\newblock Integers represented as the sum of one prime, two squares of primes
  and powers of 2.
\newblock {\em Proc. Amer. Math. Soc.}, 137(4):1185--1191, 2009.

\bibitem{PintzRom}
J.~Pintz.
\newblock A note on {R}omanov's constant.
\newblock {\em Acta Math. Hungar.}, 112(1-2):1--14, 2006.

\bibitem{PintzR}
J.~Pintz and I.~Z. Ruzsa.
\newblock On {L}innik's approximation to {G}oldbach's problem. {I}.
\newblock {\em Acta Arith.}, 109(2):169--194, 2003.

\bibitem{Ren03}
X.~M. Ren.
\newblock Sums of four cubes of primes.
\newblock {\em J. Number Theory}, 98:156--171, 2003.

\bibitem{Revol2005}
N.~Revol and F.~Rouillier.
\newblock {Motivations for an arbitrary precision interval arithmetic and the
  MPFI library}.
\newblock {\em Reliab. Comput.}, 11(4):275--290, 2005.

\bibitem{Donghua}
D.-H. Wu.
\newblock An improvement of {J}. {R}. {C}hen's theorem.
\newblock {\em Shanghai Keji Daxue Xuebao}, (1):94--99, 1987.

\bibitem{Wu}
J.~Wu.
\newblock {C}hen's double sieve, {G}oldbach's conjecture and the twin prime
  problem.
\newblock {\em Acta Arith.}, 114(3):215--273, 2004.

\bibitem{Zhao13}
L.~Zhao.
\newblock Four squares of primes and powers of 2.
\newblock {\em Acta Arith.}, 162(3):255--271, 2014.

\end{thebibliography}

\end{document}